\documentclass[10pt, a4paper]{amsart}

\usepackage{mathrsfs}
\usepackage{amsmath}
\usepackage{amssymb}
\usepackage[all]{xy}
\xyoption{all}

\newcommand{\m}{\Lambda}

\newcommand{\s}{\hfill\blacksquare}

\newcommand{\E}{\operatorname{E}}
\newcommand{\Eev}{\operatorname{E}^{\operatorname{ev}}}
\newcommand{\Eod}{\operatorname{E}^{\operatorname{odd}}}
\newcommand{\Ext}{\operatorname{Ext}\nolimits}

\newcommand{\Gr}{\operatorname{Gr}}

\newcommand{\gr}{\operatorname{gr}}

\newcommand{\Hom}{\operatorname{Hom}}

\def \K  {\text{Koszul}~}
\def \dK {~$d$\text{-Koszul}~}

\newtheorem{theorem}{Theorem}[section]

\newtheorem{definition}[theorem]{Definition}
\newtheorem{example}[theorem]{Example}

\newtheorem{lemma}[theorem]{Lemma}

\newtheorem{proposition}[theorem]{Proposition}

\newtheorem{remark}[theorem]{Remark}

\begin{document}

\title{Generalized $d$-Koszul modules}
\author [N.Bian, Y.Ye, P.Zhang]{Ning Bian, Yu Ye, Pu Zhang$^*$}

\thanks {{\it Key words and phases.} $d$-Koszul algebra, $d$-Koszul module,
generalized $d$-Koszul module.}
\thanks {$^*$The corresponding author.}
\thanks {Supported by the NSFC (10725104 and 10971206) and STCSM (09XD1402500).}
\thanks{mathbn$\symbol{64}$hotmail.com, \ yeyu$\symbol{64}$ustc.edu.cn, \ pzhang$\symbol{64}$sjtu.edu.cn}
\maketitle
\dedicatory{}%
\commby{}%
\maketitle
\begin{abstract}
Generalized \dK modules are introduced to solve an open problem: the
odd Ext-module $\Eod(M)$  of a \dK module $M$ over a $d$-Koszul
algebra $\m$ is a \K module over the even Yoneda algebra $\Eev(\m)$.
\end{abstract}

\section*{\bf Introduction}

For an integer $d\ge 2$, a \dK algebra was introduced and studied by
R. Berger \cite{B1}, and developed by E. L. Green et al. \cite{GMMZ}
to the nonlocal case. If $d=2$ it is the usual Koszul algebra. This
class of generalized Koszul structures turns out to be important for
example in theory of the Artin-Shelter algebras, the Calabi-Yau
algebras, and the Yang-Mills algebras (see e.g. [B1], [B2], [CD]).

\vskip5pt

Let $\m$ be a \dK algebra and $M$ a \dK $\m$-module. It was shown in
Theorem 6.1 of \cite{GMMZ} that the even Ext-algebra $\Eev(\m)$ is a
Koszul algebra and the even Ext-module $\Eev(M)$ is a Koszul
$\Eev(\m)$-module. This generalizes the corresponding result of J.
Backelin and R. Fr\"{o}berg \cite{BF} on the local Koszul algebras.
An open problem was raised by E. L. Green et al. \cite{GMMZ},
Section 6: Is the odd Ext-module $\Eod(M)$ also a Koszul module over
$\Eev(\m)$? E. N. Marcos and R. Mart\'inez-Villa \cite{MM} proved
that this is the case if the orthogonal algebra $\m^{!}$ is also a
\dK algebra. However, in general $\m^{!}$ is not a \dK algebra (see
\cite{B1}; also Example 2 in \cite{MM}). So the problem remains to
be open.

\vskip5pt

In this paper we introduce the so-called generalized \dK modules.
This is a natural class of graded modules. For example, the syzygies
of a \dK module are generalized \dK modules up to shifts. Also for
each $i$, $J^iM$ is a generalized \dK module up to shift, where $M$
is a generalized \dK module over a \dK algebra, and $J$ is the
graded Jacobson radical.

\vskip5pt

Our main result is as follows.

\vskip5pt

\noindent {\bf Main Theorem.} {\it Let $\m$ be a \dK algebra, and
$M$ a generalized \dK $\m$-module. Then $\Eev (M)$ is a \K module
over the \K algebra $\Eev(\m)$.}

\vskip5pt As a consequence, we have

\vskip5pt

\noindent {\bf Corollary.} {\it Let $\m$ be a \dK algebra, and $M$ a
\dK $\m$-module.  Then $\Eod (M)$ is a \K module over the \K algebra
$\Eev(\m)$.}

\vskip5pt

This answers in the affirmative the open problem mentioned above.

\section{\bf Preliminaries}

We fix the notations and recall some facts frequently used later.
For the details we refer to \cite{BGS}, \cite{GM}, and \cite{GMMZ}.

\vskip5pt

\subsection{} Throughout $\m$ is a {\it standardly graded algebra} over
a field $k$ (see [GM], p.250), i.e., $\m=\m_0\oplus \m_1\oplus
\cdots$ is a positively graded $k$-algebra satisfying the following
three conditions:

\vskip5pt

$(i)$ \ $\Lambda_0=k^r$ for some integer $r\ge 1$,

\vskip5pt

$(ii)$ \ $\dim_k \m_i< \infty, \ \forall \ i\ge 0$,

\vskip5pt

$(iii)$ \ $\m_i \m_j=\m_{i+j}, \ \forall \ i, j\ge 0$.

\vskip5pt

A left {\it graded $\m$-module} $M$ is a $\m$-module together with a
decomposition of $k$-spaces $M=\bigoplus\limits_{i\in \mathbb Z}
M_i$ such that $\Lambda_i M_j\subseteq M_{i+j}, \ \forall \ i,
j\in\mathbb Z$. Let $M$ and $N$ be graded $\m$-modules.
 A $\Lambda$-homomorphism $f:  M\to N$ is {\it a graded homomorphism}  if $f(M_i)\subseteq N_{i}, \ \forall
 \ i\in \mathbb{Z}$. For $M\in\Gr(\m)$, let
$M[n]$ denote the graded module with $M[n]_i=M_{i-n}$. Let $\m$-Mod
be the category of the left $\m$-modules,  $\Gr(\m)$ the category of
the left graded $\m$-modules and graded homomorphisms, and $\gr(\m)$
the full subcategory of $\Gr(\m)$ consisting of finitely generated
$\m$-modules. Then $\m$-Mod and $\Gr(\m)$ are abelian categories;
and $\gr(\m)$ is abelian if $\m$ is noetherian. Let $\Hom_{\Gr(\m)}$
and $\Ext^i_{\Gr(\m)}$ denote the homomorphisms and extensions in
$\Gr(\m)$, as opposed to the usual $\Hom_\m$ and $\Ext^i_\m$ in
$\m$-Mod.

\vskip5pt

Let $I$ be a subset of $\Bbb Z$, and $M\in\Gr(\m)$. $M$ is {\it
generated in degrees in $I$}, if $M = \m (\bigoplus\limits_{j\in I}
M_j)$; $M$ is {\it generated in degree $i$} if $M = \m M_i$; $M$ is
{\it supported above degree} $n$ if $M_j=0$ for $j<n$; and $M$ is
{\it concentrated in degrees in $I$} if $M_i = 0$ for $i\notin I$.

\vskip5pt

Let $J$ be the ideal $\bigoplus\limits_{i\geqslant 1}\m_i$ of $\m$.
The trivial $\m$-module $\m_0$ is the lift of the $\m_0$-module
$\m_0$ via the $k$-algebra homomorphism $\m \rightarrow \m/J=\m_0$.
It is a graded $\m$-module concentrated in degree $0$. We need the
following well-known fact.

\vskip5pt

\begin{lemma}  Let $M\in\Gr(\m)$, and $I$ be a subset of $\Bbb Z$.
If $\Hom_{\Gr(\m)}(M, \m_0[j])=0$ for $j\notin I$, then $M$ is
generated in degrees in $I$. \end{lemma} \noindent{\bf Proof.} \ For
the convenience of the reader we include a justification. Put $L: =
M/\m (\bigoplus\limits_{i\in I}M_i)$. If $L\ne 0$, then $L/JL\ne 0$.
While $L/JL$ is a graded module over the semisimple algebra $\m_0$,
it follows that $\Hom_{\Gr(\m)}(L/JL, \m_0[j])\ne 0$ for some
$j\notin I$, and hence $\Hom_{\Gr(\m)}(M, \m_0[j])\ne 0$ for some
$j\notin I$, contrary to the assumption. $\s$

\vskip5pt

\subsection{} Denote by $\E(\m)$ the Ext-algebra
$\bigoplus\limits_{i\geqslant 0}\Ext_{\m}^i(\m_0, \m_0)$, with the
multiplication given by the Yoneda product. We also consider the
even Ext-algebra $\Eev(\m): = \bigoplus\limits_{i\geqslant
0}\Ext_{\m}^{2i}(\m_0, \m_0)$, which is a positively graded algebra
with grading $\Eev(\m)_n:=\Ext^{2n}_{\m}(\m_0, \m_0)$. For a
$\m$-module $M$, let $\E(M)$ be the graded $\E({\m})$-module
$\bigoplus\limits_{i\geqslant 0}\Ext_{\m}^i(M, \m_0)$. We also
consider the even Ext-module $\Eev(M):=\bigoplus\limits_{n\geqslant
0}\Ext^{2n}_{\m}(M, \m_0)$ over $\E(\m)$, and the odd Ext-module
$\Eod(M):=\bigoplus\limits_{n\geqslant 0} \Ext^{2n+1}_{\m}(M, \m_0)$
over $\Eev(\m)$: they are graded modules with gradings
$$\Eev(M)_n:=\Ext^{2n}_{\m}(M, \m_0), \ \ \ \mbox{and} \ \ \ \Eod(M)_n:=
\Ext^{2n+1}_{\m}(M, \m_0), \ \forall \ n\ge 0.$$

\vskip5pt

Every graded $\m$-module $M$ has a graded projective resolution
\begin{align}\mathbf{Q}^{\bullet}:  \ \cdots \to Q^i \to \cdots \to Q^1 \to
Q^0 \to M \to 0.\end{align} If each $Q^i$ is finitely generated,
then we say that $\mathbf{Q}^{\bullet}$ is  {\it a finitely
generated graded projective resolution} of $M$. If $M\in\gr(\m)$,
then $M$ admits a {\it minimal} graded projective resolution $(1)$
in the sense that ${\rm Im} (Q^i\to Q^{i-1}) \subseteq JQ^{i-1}, \
\forall \ i\ge 1$ (see Propositions 2.3 and 2.4 in [GM]).

\vskip5pt

If $M\in\gr(\m)$, then for each $N\in\Gr(\m)$, $\Hom_{\m}(M, N)$ is
a graded $k$-space with {\it the shift-grading}: $\Hom_{\m}(M, N)_i
= \Hom_{\Gr(\m)}(M, N[i]), $ i.e., $\Hom_{\m}(M, N) =
\bigoplus\limits_{i\in\Bbb Z}\Hom_{\Gr(\m)}(M, N[i]).$

If $M$ has a finitely generated graded projective resolution, then
for each $N\in\Gr(\m)$ and each $n\ge 1$, $\Ext^n_{\m}(M, N)$ is a
graded $k$-space with {\it the shift grading}: $\Ext^n_{\m}(M, N)_i
=\Ext^n_{\Gr(\m)}(M, N[i]),$ i.e., $\Ext^n_{\m}(M, N) =
\bigoplus\limits_{i\in\Bbb Z}\Ext^n_{\Gr(\m)}(M, N[i]).$

\vskip5pt

Fix a minimal graded projective resolution of the trivial
$\m$-module $\m_0$:
\begin{align}\mathbf{P}^{\bullet}: \cdots \to P^n \to \cdots \to P^1 \to P^0 \to \m_0 \to 0. \end {align}

\vskip5pt

We need the following fact.

\begin{lemma}\label{supp-above} {\rm (\cite{GMMZ}, Lemma 3.2)}
Let  $M$ be a  graded module supported above degree $0$ with a
minimal graded projective resolution $(1)$. For any integer
$n\geq1$, if $P^n$ in $(2)$ is supported above degree $s$, then so
is $Q^n$.
\end{lemma}

\vskip5pt

\subsection{} For the theory of the Koszul algebras and the Koszul
modules we refer to A. Beilinson, V. Ginzburg and W. Soergel [BGS],
and E. L. Green and R. Mart\'inez-Villa [GM].

\vskip5pt

\begin{definition} [\cite{GMMZ}, \cite{MM}] Let $d\ge 2$ be an integer.  A graded $\m$-module $M$
is a $d$-Koszul module if $M$ admits a finitely generated graded
projective resolution $(1)$ such that each $Q^i$ is generated in
degree $\delta(i)$, where
$$\delta(i): =
\begin{cases} nd, & \ \mbox{if~}i=2n,\\
nd+1, & \ \mbox{if~}i=2n+1.
\end{cases} $$

\vskip5pt

If the trivial $\m$-module $\m_0$ is a \dK module, then we call $\m$
a \dK algebra.
\end{definition}

\vskip10pt

\begin{theorem} \label{GMMZthm} {\rm (\cite{GMMZ}, Theorem 6.1)}
Let $\m$ be a \dK algebra and $M$ a \dK $\m$-module. Then $\Eev(\m)$
is a \K algebra, and $\Eev(M)$ is a \K $\Eev(\m)$-module.
\end{theorem}

\section{\bf Generalized $d$-Koszul modules}

\subsection{} Let $d\ge 2$ be an integer. For each integer $i\ge 0$ we
assign a subset $\Delta(i)$ of $\Bbb N_0$ as
$$\Delta(i): =
\begin{cases} \{nd\}, & \ \mbox{if}~i=2n;\\
\{nd+1, \cdots, nd+d-1\}, & \ \mbox{if}~i=2n+1.
\end{cases}$$

\vskip10pt

\begin{definition} A graded $\m$-module $M$ is called a {\it generalized \dK module} if $M$ admits a
finitely generated graded projective resolution
$\mathbf{Q}^{\bullet}$  such that each $Q^i$ is generated in degrees
in $\Delta(i)$, i.e., $Q^i=\m
(\bigoplus\limits_{j\in\Delta(i)}Q_j^i), \ i\ge 0$.
\end{definition}

\vskip5pt

\begin{remark} $(i)$ \ As remarked by Beilinson-Ginzburg-Soergel [BGS] (p.476)
in the Koszul situation, $\mathbf{Q}^{\bullet}$ in Definition 2.1 is
unique up to isomorphism.

More precisely, if $\mathbf{L}^{\bullet}$ is another graded
projective resolution of $M$ such that each $L^i$ is also generated
in degrees in $\Delta(i)$ (it is not assumed to be finitely
generated), then $\mathbf{L}^{\bullet}\cong \mathbf{Q}^{\bullet}$ as
complexes. In fact, $\mathbf{L}^{\bullet}$ is homotopy equivalent to
$\mathbf{Q}^{\bullet}$; while any chain maps
$\mathbf{Q}^{\bullet}\rightarrow \mathbf{Q}^{\bullet}$ and
$\mathbf{L}^{\bullet}\rightarrow \mathbf{L}^{\bullet}$, which
respect the grading on $Q^i$ and $L^i$ and are homotopic to zero
must themselves be zero (since any element in $\Delta(i)$ is
strictly smaller than any element in $\Delta(i+1)$, and $Q^i$ and
$L^i$ are both generated in degrees in $\Delta(i)$). It follows that
$\mathbf{L}^{\bullet}\cong \mathbf{Q}^{\bullet}$ as complexes.

\vskip5pt

 $(ii)$ \ We emphasize that, as in
the $d$-Koszul situation, here $Q^i$ is also required to be {\bf
finitely generated}: it is for the application of the shift grading
on $\Ext^n_\m(M, -)$.

\vskip5pt

$(iii)$ \ If $M$ is a generalized $d$-Koszul module, then such a
graded projective resolution $\mathbf Q^\bullet$ in the definition
is minimal, and each syzygy $\Omega^i(M)$ is a graded $\m$-module
finitely generated in degrees in $\Delta(i)$. In particular, $M$ is
finitely generated in degree $0$.

\vskip5pt

$(iv)$ \ A \dK module is always generalized $d$-Koszul; and a
generalized $2$-Koszul module is a finitely generated Koszul module
(if $\m$ is noetherian, then a generalized $2$-Koszul $\m$-module is
exactly a finitely generated Koszul $\m$-module).
\end{remark}

\vskip5pt

\begin{example} \ Let $A$ be the algebra given by
the quiver
\[
\xymatrix{1 \ar@(dl, ul)^\alpha \ar[r]^\beta & 2 \ar[r]^\gamma & 3 }
\]
with relations $\alpha^3, \ \gamma \beta \alpha$. Then the simple
(left) module $S(1)$ has a minimal graded projective resolution
$$\cdots \rightarrow P(1)[4]\oplus P(2)[5]\rightarrow P(1)[3]\oplus P(3)[3]\rightarrow P(1)[1]\oplus P(2)[1]
\rightarrow P(1)\rightarrow S(1)\rightarrow 0$$ where $\Omega^4 S(1)
= (\Omega^2S(1))[3]$. Thus $S(1)$ is a generalized $3$-\K
$A$-module. Since $Q^3 = P(1)[4]\oplus P(2)[5]$ is generated in
degrees $4$ and $5$, but {\it not} generated in degree $4$, it
follows that $S(3)$ is {\bf not} a $3$-Koszul $A$-module (by an
argument in Remark 2.2$(i)$).
\end{example}

\vskip5pt

\subsection{} We have the following characterization for a $d$-Koszul
module and for a generalized $d$-Koszul module, which is the
corresponding version of Proposition 2.14.2 in Beilinson - Ginzburg
- Soergel [BGS] for the Koszul modules.

\vskip5pt
\begin{lemma}\label{gdk char} Let $M$ be a graded $\m$-module with
a finitely generated graded projective resolution. Then

\vskip5pt

$(i)$ \ $M$ is \dK if and only if $\Ext_{\m}^i (M, \m_0)_j =
\Ext_{\Gr(\m)}^i (M, \m_0[j])=0, \ \forall \ j\ne
\delta(i)$.

\vskip5pt

$(ii)$ \ $M$ is generalized \dK if and only if $\Ext_{\m}^i (M,
\m_0)$ is concentrated in degrees in $\Delta(i)$, with the shift
grading, i.e., \ $\Ext_{\m}^i (M, \m_0)_j = \Ext_{\Gr(\m)}^i (M,
\m_0[j])=0, \ \forall \ j\notin \Delta(i).$
\end{lemma}
\noindent{\bf Proof.} \ They can be similarly proved as Proposition
2.14.2 in [BGS]. For the convenience of the reader we include a
justification of $(ii)$.

Assume that $M$ is generalized $d$-Koszul. Then $M$ has a graded
projective resolution $\mathbf Q^\bullet$ such that each $Q^i$ is
generated in degrees in $\Delta(i)$, and $\Ext_{\Gr(\m)}^i (M,
\m_0[j])$ is the $i$-th cohomology group of the complex
$\Hom_{\Gr(\m)}(\mathbf{Q}^{\bullet}, \m_0[j])$. Since $Q^i$ is
generated in degrees in $\Delta(i)$, and $\m_0[j]$ is concentrated
in degree $j$, it follows that $\Hom_{\Gr(\m)}(Q^i, \m_0[j])=0$ for
$j\notin\Delta(i)$, and hence $\Ext_{\Gr(\m)}^i (M, \m_0[j])=0$ for
$j\notin \Delta(i)$.

\vskip5pt

Conversely, assume that $\Ext_{\Gr(\m)}^i (M, \m_0[j])=0$ for
$j\notin \Delta(i)$. We construct inductively a graded projective
resolution $\mathbf{L}^{\bullet}$ of $M$ such that each $L^i$ is
generated in degrees in $\Delta(i)$. Since $\Hom_{\Gr(\m)}(M,
\m_0[j])=0$ for $j\ne 0$, by Lemma 1.1, $M$ is generated in degree
$0$, and hence we have a surjective graded $\m$-homomorphism
$L^0\rightarrow M$ such that $L^0$ is generated in degree $0$.
Denote by $K^1$ its kernel. Then $\Hom_{\Gr(\m)}(K^1,
\m_0[j])=\Ext_{\Gr(\m)}^1(M, \m_0[j]) =0$ for $j\notin \Delta(1)$,
and hence by Lemma 1.1, $K^1$ is generated in degrees in
$\Delta(1)$. Thus we have a surjective graded $\m$-homomorphism
$L^1\rightarrow K^1$ such that $L^1$ is generated in degrees in
$\Delta(1)$. Repeating this process we are done.

\vskip5pt

By assumption we have already a finitely generated graded projective
resolution $\mathbf{Q}^{\bullet}$. By the argument in Remark
2.2$(i)$, there are chain maps $f: \mathbf{L}^{\bullet}\rightarrow
\mathbf{Q}^{\bullet}$ and $g: \mathbf{Q}^{\bullet}\rightarrow
\mathbf{L}^{\bullet}$ such that $gf = {\rm
Id}_{\mathbf{L}^{\bullet}}$, which means that $\mathbf{L}^{\bullet}$
is a direct summand of $\mathbf{Q}^{\bullet}$. Thus
$\mathbf{L}^{\bullet}$ is also a finitely generated resolution. By
definition $M$ is generalized $d$-Koszul. $\s$

\vskip15pt

\subsection{} For a $d$-Koszul module $M$, in general $\Omega^iM$ and
$J^iM$ are {\bf not} \dK modules, up to shifts (see Proposition 5.2
in [GMMZ] for some special cases); however, they turn out to be
generalized $d$-Koszul, after proper shifts. In the rest of this
section we precisely state and prove these results, which will be
important in the proof of Main Theorem and Corollary.

\begin{lemma} \label{lemma1} \ Let $M$ be a \dK $\m$-module. Then

\vskip5pt

$(i)$ \  $\Eod(M)\cong\Eev(\Omega M)$ as graded $\E(\m)$-modules.

\vskip5pt

$(ii)$ \  $(\Omega^i M)[-\delta(i)]$ is a generalized \dK module for
each $i\ge 0$.
\end{lemma}
\noindent{\bf Proof.} $(i)$ \ By definition we have an isomorphism
of graded $\E(\m)$-modules
$$\Eod(M)=\bigoplus_{n\geqslant 0} \Ext^{2n+1}_{\m}(M, \m_0)
\cong \bigoplus_{n\geqslant 0} \Ext^{2n}_{\m}(\Omega M, \m_0) =
\Eev(\Omega M).$$

$(ii)$ \  Taking a graded projective resolution
$\mathbf{Q}^{\bullet}$ of $M$ such that each $Q^i$ is finitely
generated in degrees in $\delta(i)$, we see that $(\Omega^i
M)[-\delta(i)]$ has a graded projective resolution
$$\mathbf{L}^{\bullet}:  \ \cdots \to L^j \to \cdots \to L^1 \to L^0
\to (\Omega^i M)[-\delta(i)] \to 0$$ where $L^j =
Q^{i+j}[-\delta(i)]$ is finitely generated in degree $\delta(i+j) -
\delta(i)$ for $j\ge 0$.

If  $i$ is even, then $L^j$ is generated in degree $\delta(j)$. That
is, $(\Omega^i M)[-\delta(i)]$ is a \dK module, and hence a
generalized \dK module.

Assume that $i$ is odd. Then $L^j$ is generated in degree $nd$ if
$j= 2n$,  and $L^j$ is generated in degree $nd + d -1$ if $j= 2n+1$.
By definition $(\Omega^i M)[-\delta(i)]$ is generalized \dK. $\s$

\begin{theorem} \label{radical-gdk}
Let  $\m$ be a  \dK algebra and $M$ a generalized \dK $\m$-module.
Then \vskip5pt

$(i)$ \  $(J^iM)[-i]$ is generalized \dK for each $i\ge 1$.
\vskip5pt

$(ii)$ \ For each $n\ge 1$ we have $k$-isomorphisms
$$\Ext^{2n-1}_{\Gr(\m)}(JM, \m_0[nd])\cong\Ext^{2n-1}_{\Gr(\m)}(J^2M, \m_0[nd])\cong
\cdots \cong\Ext^{2n-1}_{\Gr(\m)}(J^{d-1} M, \m_0[nd]).$$
\end{theorem}

\noindent{\bf Proof.} $(i)$ \ It suffices to prove that $(JM)[-1]$
is generalized $d$-Koszul. Since $M = \bigoplus\limits_{i\geqslant
0}M_i$ is finitely generated in degree $0$,
$JM=\bigoplus\limits_{i\geqslant 1}M_i$ is finitely generated in
degree $1$.

\vskip5pt We first prove the following claim: $JM[-1]$ admits a
graded projective resolution $\mathbf Q_1^\bullet$ such that $Q_1^i$
is generated in degrees in $\Delta(i)$. By the proof of Lemma
2.4$(ii)$, it suffices for each $n\ge 0$ to prove that
$\Ext^{2n}_{\Gr(\m)}(JM[-1], \m_0[j]) = \Ext^{2n}_{\Gr(\m)}(JM,
\m_0[j+1]) = 0, \ \ \forall \ j\ne nd,$ and that
$\Ext^{2n+1}_{\Gr(\m)}(JM, \m_0[j+1])=0, \ \ \forall \  j\notin
\Delta(2n+1) = \{nd+1, \cdots, nd+d-1\}.$

\vskip5pt

Applying $\Hom_{\Gr(\m)}(-, \m_0[j+1])$ to the graded  exact
sequence $0 \to JM \to M \to M/JM \to 0$ we get the following exact
sequence of $k$-spaces {\small{\begin{align*} \Ext^{2n}_{\Gr(\m)}(M,
\m_0[j+1])\to \Ext^{2n}_{\Gr(\m)}(JM, \m_0[j+1])\to
 \Ext^{2n+1}_{\Gr(\m)}(M/JM, \m_0[j+1])
 \\ \to \Ext^{2n+1}_{\Gr(\m)}(M, \m_0[j+1])\to
\Ext^{2n+1}_{\Gr(\m)}(JM, \m_0[j+1])\to
\Ext^{2(n+1)}_{\Gr(\m)}(M/JM, \m_0[j+1]).\end{align*}}}\hskip8pt
Since $\m$ is a $d$-Koszul algebra, $P^{2n}$ in $(2)$ is supported
above degrees $nd$, and hence by Lemma \ref{supp-above}, $Q^{2n}$ is
supported above degrees $nd$, where $\mathbf Q^\bullet$ is a minimal
graded projective resolution of $JM[-1]$. Thus
$\Ext^{2n}_{\Gr(\m)}(JM, \m_0[j+1]) = 0$ for $j < nd$. Similarly,
$\Ext^{2n+1}_{\Gr(\m)}(JM, \m_0[j+1]) = 0$ for $j < nd+1$.

\vskip5pt

Since $M$ is generalized $d$-Koszul, by Lemma 2.4$(ii)$,
$\Ext^{2n}_{\Gr(\m)}(M, \m_0[j+1])=0$ if $j\ne nd-1$,  and
$\Ext^{2n+1}_{\Gr(\m)}(M, \m_0[j+1]) =0$ if $j\notin \{nd, \cdots,
nd+d-2\}$.

\vskip5pt

Note that $M/JM$ is a $\m/J$-module and $\m/J = \m_0$ is a
semisimple algebra. Thus $M/JM$ is a direct summand of a finite
direct sum of copies of the trivial $\m$-module $\m_0$. In
particular, $M/JM$ is a $d$-Koszul module. By Lemma 2.4$(i)$,
$\Ext^{2n+1}_{\Gr(\m)}(M/JM, \m_0[j+1]) =0 \ \ \mbox{if} \ j\ne nd,$
and $\Ext^{2(n+1)}_{\Gr(\m)}(M/JM, \m_0[j+1])=0 \ \ \mbox{if} \ j\ne
(n+1)d-1.$

\vskip5pt

Now if $j\ne nd$, then by the exact sequence above we have the
following exact sequence {\small
{\begin{align*}\Ext^{2n}_{\Gr(\m)}(M, \m_0[j+1])\to
\Ext^{2n}_{\Gr(\m)}(JM, \m_0[j+1])\to \Ext^{2n+1}_{\Gr(\m)}(M/JM,
\m_0[j+1])=0,\end{align*}}}\hskip-1pt where if  $j\ne nd-1$ then
$\Ext^{2n}_{\Gr(\m)}(M, \m_0[j+1]) = 0$, and hence
$\Ext^{2n}_{\Gr(\m)}(JM, \m_0[j+1])=0$; and if $j=nd-1 < nd$, then
we already know $\Ext^{2n}_{\Gr(\m)}(JM, \m_0[j+1])=0$.

\vskip5pt

Let $j\notin \Delta(2n+1) = \{nd+1, \cdots, nd+d-1\}.$ Then by the
exact sequence above we have the following exact sequence {\small
{\begin{align*}\Ext^{2n+1}_{\Gr(\m)}(M, \m_0[j+1])\to
\Ext^{2n+1}_{\Gr(\m)}(JM, \m_0[j+1])\to
\Ext^{2(n+1)}_{\Gr(\m)}(M/JM, \m_0[j+1])=0,\end{align*}}} \hskip-1pt
where if $j\notin \{nd, \cdots, nd+d-2\}$ then
$\Ext^{2n+1}_{\Gr(\m)}(M, \m_0[j+1]) = 0$, and hence
$\Ext^{2n+1}_{\Gr(\m)}(JM, \m_0[j+1])=0$; and if $j\in \{nd, \cdots,
nd+d-2\}$, then $j = nd < nd+1$, and in this case we already know
$\Ext^{2n+1}_{\Gr(\m)}(JM, \m_0[j+1])=0$.  This proves the claim.

\vskip5pt

Since $M/JM$ is a $d$-Koszul module, $M/JM$ has a finitely generated
graded projective resolution, say $\mathbf Q_3^\bullet$, such that
$\mathbf{Q}_3^{i}$ is generated in degrees in $\delta(i)$. By the
graded version of the Horseshoe Lemma, we get a graded projective
resolution $\mathbf Q_2^\bullet$ of $M$, such that
$\mathbf{Q}_2^{i}= \mathbf{Q}_1^{i}\oplus \mathbf{Q}_3^{i}$ for each
$i$. Thus  $\mathbf{Q}_2^{i}$ is also generated in degrees in
$\Delta(i)$.  Since $M$ is a generalized $d$-Koszul module, by
Remark 2.2$(i)$, we know that $\mathbf Q_2^\bullet$ is finitely
generated, and hence $\mathbf Q_1^\bullet$ is finitely generated. By
definition $JM[-1]$ is generalized $d$-Koszul.

\vskip5pt

$(ii)$ \ Let $d\ge 3$. Applying $\Hom_{\m}(-, \m_0)$ to the graded
exact sequence $0\rightarrow J^2M \rightarrow JM \rightarrow JM/J^2M
\rightarrow 0,$ we get the following exact sequence {\small
\begin{align*} \Ext^{2n-1}_{\m}(JM/J^2M, \m_0) \to
\Ext^{2n-1}_{\m}(JM, \m_0) \to \Ext^{2n-1}_{\m}(J^2M, \m_0) \to
\Ext^{2n}_{\m}(JM/J^2M, \m_0).\end{align*}} \hskip-3pt Since
$(JM/J^2M)[-1]$ is $d$-Koszul, by Lemma 2.4$(i)$,
$\Ext^{2n-1}_{\Gr(\m)}(JM/J^2M, \m_0[j]) = 0$ if $j\ne nd-d+2$,  and
$\Ext^{2n}_{\m}(JM/J^2M, \m_0[j])=0$ if $j\ne nd+1$. Taking the
$nd$-th homogeneous components of the exact sequence above, we
obtain that $\Ext^{2n-1}_{\Gr(\m)}(JM,
\m_0[nd])\cong\Ext^{2n-1}_{\Gr(\m)}(J^2M, \m_0[nd])$. Repeating the
process one gets $(ii)$.  $\s$

\section{\bf Proofs of Main Theorem and Corollary}

\subsection{} We begin with a lemma, which seems to be of independent interest.

\begin{lemma}\label{Koszul class}
Let $A$ be an arbitrary Koszul algebra and $\mathcal{C}$  a full
subcategory of \ $\Gr(A)$. Suppose that for any $X\in\mathcal{C}$,
there exist exact sequences in $\Gr(A)$
\begin{align}
0 \to \Omega \to P^0 \to X \to 0,\end{align} \begin{align}0 \to
X''\to X' \to \Omega[-1] \to 0,
\end{align}
such that  $P^0$ is a graded projective $A$-module generated in
degree $0$ and $X', X''\in \mathcal{C}$. Then all modules in
$\mathcal{C}$ are Koszul $A$-modules.

\end{lemma}

\noindent {\bf Proof.} \ By Proposition 2.14.2 in Beilinson -
Ginzburg - Soergel [BGS], it suffices to prove that for each $X\in
\mathcal{C}$, $\Ext^i_{\Gr(A)}(X, A_0[j]) = 0$ unless $j = i$. We
use induction on $i$.

The sequence $(3)$ implies that $X$ is generated in degree $0$, and
hence $\Hom_{\Gr(A)}(X, A_0[j]) = 0$ unless $j = 0$. The sequence
$(4)$ implies that $\Omega$ is a graded $A$-module and is generated
in degree $1$, since $X'\in\mathcal C$ is generated in degree $0$.
By $\Ext^1_{\Gr(A)}(X, A_0[j])\cong \Hom_{\Gr(A)}(\Omega, A_0[j])$,
we see that $\Ext^1_{\Gr(A)}(X, A_0[j]) = 0$ unless $j =1$.

Let $n\ge 1$. Assume that for each $X\in \mathcal{C}$ and for each
positive integer $i$ with $i\le n$, $\Ext^i_{\Gr(A)}(X, A_0[j]) =0$
unless $j = i$. The exact sequence $(4)$ implies the following exact
sequence for every integer $j$
$$\Ext^{n-1}_{\Gr(A)}(X'', A_0[j])\to \Ext^n_{\Gr(A)}(\Omega[-1], A_0[j]) \to \Ext^{n}_{\Gr(A)}(X', A_0[j]).$$
By the inductive hypothesis, we have $\Ext^{n-1}_{\Gr(A)}(X'',
A_0[j]) = 0$ unless $j = n-1$, and $\Ext^{n}_{\Gr(A)}(X', A_0[j]) =
0$ unless $j = n$. Let $\mathbf Q^\bullet$ be a minimal graded
projective resolution of $\Omega[-1]$ (it exists since $\Omega
\subseteq P^0$ is supported above $0$). By Lemma \ref{supp-above},
$Q^n$ is supported above degree $n$, which implies
$\Ext^{n}_{\Gr(A)}(\Omega[-1], A_0[j])$ $ = 0$ for $j < n$. It
follows from the exact sequence above that
$\Ext^n_{\Gr(A)}(\Omega[-1], A_0[j]) = 0$ unless $j = n$. Thus
$\Ext^{n+1}_{\Gr(A)}(X, A_0[j])=\Ext^{n}_{\Gr(A)}(\Omega, A_0[j])=
\Ext^{n}_{\Gr(A)}(\Omega[-1], A_0[j-1]) = 0$ unless $j = n+1$. This
completes the proof. $\s$

\vskip10pt

\subsection{\bf Proof of Main Theorem} \ By Theorem
\ref{GMMZthm}, $\Eev(\m)$ is a Koszul algebra. Put
$$\mathcal{C}:=\{\Eev(N) \in {\Gr(\Eev(\m))} \mid\, N \ \mbox{is a generalized \dK $\m$-module} \}.$$
It suffices to prove that all the conditions in Lemma \ref{Koszul
class} are satisfied.

\vskip5pt

The graded exact sequence $0 \to JN \to N \to N/JN \to 0$ induces
the following exact sequence of graded $k$-spaces for each $n\ge 0$
{\small
\begin{align} \Ext^{2n-1}_{\m}(N, \m_0) \to \Ext^{2n-1}_{\m}(JN,
\m_0) \to \Ext^{2n}_{\m} (N/JN, \m_0) \to \Ext^{2n}_{\m}(N, \m_0)
\to \Ext^{2n}_{\m}(JN, \m_0).
\end{align}}
\hskip-3ptSince $N$ and $JN[-1]$ are generalized $d$-Koszul, by
Lemma 2.4$(ii)$, we have {\small $\Ext^{2n-1}_{\Gr(\m)}(N, \m_0[nd])
= 0 = \Ext^{2n}_{\Gr(\m)}(JN, \m_0[nd])$.}  Taking the $nd$-th
homogeneous components of $(5)$ we get the following exact sequence
for each $n\ge 0$
\begin{align} 0 \to \Ext^{2n-1}_{\Gr(\m)}(JN,
\m_0[nd]) \to \Ext^{2n}_{\Gr(\m)} (N/JN, \m_0[nd]) \to
\Ext^{2n}_{\Gr(\m)}(N, \m_0[nd]) \to 0.
\end{align}
Since $N$ is generalized $d$-Koszul and $N/JN$ is $d$-Koszul, by
Lemma 2.4, we have
$$\Eev(N/JN) = \bigoplus\limits_{n\ge
0}\Ext^{2n}_{\Gr(\m)} (N/JN, \m_0[nd]), \ \ \Eev(N) =
\bigoplus\limits_{n\ge 0}\Ext^{2n}_{\Gr(\m)} (N, \m_0[nd]).$$ By
taking direct sum of $(6)$, we get the following short exact
sequence in $\Gr(\Eev(\m))$:
\begin{equation}\label{mthm eqn1}0 \to \Omega
\to \Eev(N/JN) \to \Eev(N) \to 0,\end{equation} where $\Omega
:=\bigoplus\limits_{n\ge 0}\Ext^{2n-1}_{\Gr(\m)}(JN, \m_0[nd])$. In
particular,  $\Omega$ is a graded $\E^{ev}(\m)$-module with grading
$\Omega_n: = \Ext^{2n-1}_{\Gr(\m)}(JN, \m_0[nd])$.  (One can also
prove this directly as follows: since $\m$ is $d$-Koszul algebra, it
follows from Lemma 2.4$(i)$ that
\begin{align*}\Ext^{2m}_{\m}(\m_0, \m_0) \Ext^{2n-1}_{\Gr(\m)}(JN,
\m_0[nd])  & = \Ext^{2m}_{\Gr(\m)}(\m_0,
\m_0[md])\Ext^{2n-1}_{\Gr(\m)}(JN, \m_0[nd])\\ & \subseteq
\Ext^{2(n+m)-1}_{\Gr(\m)}(JN, \m_0[(n+m)d]).)\end{align*}

\vskip5pt

By Theorem 1.4, $\Eev(N/JN)$ is a Koszul $\Eev(\m)$-module, in
particular it is generated in degree $0$. Since $N/JN$ is a direct
summand of finite direct sum of copies of the trivial $\m$-module
$\m_0$, $\Eev(N/JN)$ is a projective $\Eev(\m)$-module.

\vskip5pt

Similarly, the graded exact sequence $0 \to J^dN \to J^{d-1}N \to
J^{d-1}N/J^dN \to 0$ induces the following exact sequence of graded
$k$-spaces for each $n\ge 0$ {\small
\begin{align*} \Ext^{2n}_{\m}(J^{d-1}N, \m_0) &\to \Ext^{2n}_{\m}(J^dN,
\m_0) \to \Ext^{2n+1}_{\m} (J^{d-1}N/J^dN, \m_0) \\ & \to
\Ext^{2n+1}_{\m}(J^{d-1}N, \m_0) \to \Ext^{2n+1}_{\m}(J^dN, \m_0).
\end{align*}}
\hskip-3pt Note that by Theorem 2.6$(i)$, $J^{d-1}N[-(d-1)]$ and
$J^dN[-d]$ are generalized $d$-Koszul $\m$-modules, and that
$(J^{d-1}N/J^dN)[-(d-1)]$ is a $d$-Koszul module. Taking the
$(n+1)d$-th homogeneous components, and by the same arguments we get
another exact sequence in $\Gr(\Eev(\m))$:
{\small\begin{align*}\label{mthm eqn2} 0 \to \bigoplus\limits_{n \ge
0} \Ext^{2n}_{\Gr(\m)}(J^{d}N, \m_0[(n+1)d]) &\to
\bigoplus\limits_{n \ge 0} \Ext^{2n+1}_{\Gr(\m)}(J^{d-1}N/J^{d}N,
\m_0[(n+1)d]) \\ &\to \bigoplus\limits_{n \ge 0}
\Ext^{2n+1}_{\Gr(\m)}(J^{d-1}N, \m_0[(n+1)d]) \to 0,
\end{align*}}
or equivalently, \begin{equation}\label{mthm eqn2} 0 \to \Eev(J^d N)
\to \Eod(J^{d-1}N/J^d N) \to \Omega[-1] \to 0,
\end{equation} where
\begin{align*}\Omega[-1] & = (\bigoplus\limits_{n\geq 0}\Ext^{2n-1}_{\Gr(\m)}(JN,
\m_0[nd]))[-1] = \bigoplus\limits_{n\geq 0}\Ext^{2n+1}_{\Gr(\m)}(JN,
\m_0[(n+1)d])\\ & \cong \bigoplus\limits_{n\geq
0}\Ext^{2n+1}_{\Gr(\m)}(J^{d-1}N, \m_0[(n+1)d]),
\end{align*}
where the last isomorphism follows from Theorem
\ref{radical-gdk}$(ii)$.

\vskip5pt

Since $(J^{d-1}N/J^{d}N)[-(d-1)]$ is $d$-Koszul, by Lemma
\ref{lemma1}$(ii)$,  $\Omega (J^{d-1}N/J^dN)[-d]$ is generalized
$d$-Koszul, and by Lemma \ref{lemma1}$(i)$,  we have
\begin{align*} \Eod(J^{d-1}N/J^dN) =  \Eod((J^{d-1}N/J^{d}N)[-(d-1)]) \cong
\Eev(\Omega(J^{d-1}N/J^dN)[-d]), \end{align*} from which we see
$\Eod(J^{d-1}N/J^dN)\in \mathcal{C}$.

Since $N$ is generalized $d$-Koszul, by Theorem
\ref{radical-gdk}$(i)$, $(J^dN)[-d]$ is generalized $d$-Koszul. Thus
$\Eev(J^d N)= \Eev((J^d N)[-d])\in \mathcal{C}$.

Now $(7)$ and $(8)$ shows that all the conditions in Lemma
\ref{Koszul class} are satisfied. This completes the proof. $\s$

\vskip10pt

\subsection {\bf Proof of Corollary} By Lemma
\ref{lemma1}$(ii)$, $(\Omega M)[-1]$ is a generalized $d$-Koszul
module. It follows from Main Theorem that $\Eev((\Omega M)[-1])$ is
a Koszul $\Eev(\m)$-module. Therefore by Lemma \ref{lemma1}$(i)$,
$\Eod (M)\cong \Eev(\Omega M) = \Eev((\Omega M)[-1])$ is a Koszul
$\Eev(\m)$-module. $\s$

\vskip15pt

{\bf Acknowledgements.} We would like to thank Edward L. Green for
his valuable conversations and discussions on the Koszulity and on
the problems discussed in this paper. We also thank Eduardo N.
Marcos and Roberto Mart\'inez-Villa for reading the manuscript and
giving comments.

\vskip20pt

{\small \noindent N. Bian: \ Dept. Math., \ Shanghai Jiao Tong
University, \ Shanghai 200240, P. R. China

\ \ \ \ \ \ \ \ \ \ Dept. Math., \ Shandong University of
Technology, \ Zibo 255049, P. R. China \vskip5pt

 \noindent Y. Ye: \
Dept. Math., \ University of Science and Technology of China, \
Hefei 230026, P. R. China

\vskip5pt

\noindent P. Zhang: \ Dept. Math., \ Shanghai Jiao Tong University,
\ Shanghai 200240, P. R. China}

\end{document}